\documentclass[3p,sort&compress]{elsarticle}
\usepackage[T1]{fontenc}
\usepackage{lmodern}
\usepackage{amsmath,amssymb,amsthm}
\usepackage{microtype}
\usepackage[hidelinks]{hyperref}
\newtheorem{thm}{Theorem}[section]
\newtheorem{theorem}[thm]{Theorem}
\newtheorem{lemma}[thm]{Lemma}
\newtheorem{proposition}[thm]{Proposition}
\newtheorem{corollary}[thm]{Corollary}
\theoremstyle{definition}

\newtheorem{remark}[thm]{Remark}

\newcommand{\sat}{\operatorname{sat}}
\newcommand{\KIrr}{\operatorname{K-Irr}}
\newcommand{\Fix}{\operatorname{Fix}}
\newcommand{\cl}{\operatorname{cl}}
\newcommand{\N}{\mathbb N}
\newcommand{\cA}{\mathcal A}

\hypersetup{pdftitle={Kernel retracts of co-sober spaces},
  pdfauthor={Xiaoquan Xu, Fengshuo Xu, Haitao Ji and Dong Qiu},
  pdfkeywords={co-sober space, compact saturated set, kernel retract, KC-space, anti-Hausdorff space}}

\begin{document}
\begin{frontmatter}
\title{Kernel retracts of co-sober spaces\tnoteref{funding}}
\tnotetext[funding]{This research was supported by the National Natural Science Foundation of China (Nos. 12471070, 12571489, 12071199) and Guangxi Natural Science Foundation (No. 2025GXNS-FAA069576).}
\author[nit]{Xiaoquan Xu\corref{cor1}}
\ead{xiqxu2002@163.com}
\author[gxu]{Fengshuo Xu}
\ead{18724184305@163.com}
\author[gxu]{Haitao Ji}
\ead{2411130238@st.gxu.edu.cn}
\author[gxu]{Dong Qiu}
\ead{qiudong@gxu.edu.cn}
\cortext[cor1]{Corresponding author.}
\address[nit]{School of Computer Information Engineering, Nanchang Institute of Technology, Nanchang 330044, China}
\address[gxu]{School of Mathematics, Guangxi University, Nanning 530004, China}

\begin{abstract}
Kernel retracts of co-sober spaces need not be co-sober, even when the ambient space is countably infinite and compact. We construct a co-sober $T_0$-space $Y$ and a continuous idempotent map $k\colon Y\to Y$ satisfying $k(y)\leq_Y y$ in the specialization order, whose image with the subspace topology is not co-sober. The construction starts from a countable compact $T_1$ KC anti-Hausdorff space of van Douwen type and adjoins an infinite chain. We classify the $k$-irreducible compact saturated subsets of the extension and show that each is a point saturation. The kernel operator collapses the added chain to a distinguished point while fixing the original space, thereby recovering the non-co-sober seed as its image. This answers negatively two questions of Xu and Zhao, concerning continuous kernel operators and arbitrary retracts, and contrasts with the known preservation of co-sobriety by continuous closure operators.
\end{abstract}
\end{frontmatter}

\section{Introduction}

Sobriety and co-sobriety are distinct properties in non-Hausdorff topology and domain theory. The notion of a co-sober space was introduced by Escard\'o, Lawson and Simpson in connection with a dual form of the Hofmann--Mislove theorem; see \cite{EscardoLawsonSimpson2004}. In a co-sober space, every $k$-irreducible compact saturated set is the saturation of a unique point. It is easy to give a co-sober space which is not sober (see, e.g.,  \cite[Example 2.1]{WenXu2018}). Escard\'o, Lawson and Simpson asked whether every sober space is co-sober (\cite[Problem 9.7]{EscardoLawsonSimpson2004}). In \cite{WenXu2018}, Wen and Xu gave a negative answer to this problem by proving that Isbell's lattice constructed in \cite{Isbell-1982} equipped with the lower topology is sober but not co-sober. Based on Wen and Xu's result, Shen et al. \cite{Shen-Wu-Xi-Zhao-2020} further obtained a dcpo whose Scott topology is sober but not co-sober.

A basic theme in the theory is the preservation of co-sobriety under natural topological constructions (see \cite{WenXu2018, XuZhao2021, Zhao-He-Wang-2023}). Wen and Xu proved, among other things, that closed subspaces and saturated subspaces of co-sober spaces are co-sober, and that the image of a co-sober space under a continuous closure operator is co-sober \cite{WenXu2018}. Motivated by this result, Xu and Zhao asked whether the analogous assertion holds for kernel operators \cite[Question~9.5]{XuZhao2021}. With the convention
\[
x\leq_X y ~\Longleftrightarrow~ x\in\cl_X\{y\},
\]
a continuous map $k\colon X\to X$ is a kernel operator if
\[
k^2=k,~ k(x)\leq_X x~(x\in X).
\]
They also asked whether an arbitrary retract of a co-sober space must be co-sober \cite[Question~9.6]{XuZhao2021}.

We answer both questions negatively within the class of countable compact spaces. The construction separates two requirements: co-sobriety of the ambient space and failure of co-sobriety in the kernel image. We start with a countable compact $T_1$ KC anti-Hausdorff seed $Z$ which is not co-sober. We then attach a countable chain to $Z$. The chain is topologized so that all non-principal compact saturated sets whose projection covers the whole seed admit explicit nontrivial decompositions and hence fail to be $k$-irreducible. The remaining $k$-irreducible compact saturated sets project onto proper $k$-irreducible compact subsets of $Z$, and those are forced to be singletons. This makes the extension co-sober. Finally, a kernel operator collapses the added chain to one distinguished point of $Z$ while fixing the old copy of $Z$, so that its image recovers the non-co-sober seed.

It is useful to formulate the conclusion in the language of retracts. If $r\colon X\to A$ and $i\colon A\to X$ are continuous, $r i=\mathrm{id}_A$, and $i r\leq\mathrm{id}_X$, we shall call $A$ a \emph{kernel retract} of $X$. Dually, the condition $\mathrm{id}_X\leq i r$ corresponds to the closure-operator situation; one may naturally call such an $A$ a \emph{closure retract}. The known positive theorem for closure operators and our negative theorem for kernel operators therefore exhibit a genuine asymmetry between these two kinds of order-compatible retracts.

Our main result is the following.

\begin{theorem}\label{thm:main-intro}
There exist a countably infinite compact co-sober $T_0$-space $Y$ and a continuous kernel operator $k\colon Y\to Y$ such that the subspace $k(Y)$ is not co-sober. Consequently, a kernel retract, and hence a retract, of a co-sober space need not be co-sober.
\end{theorem}

The proof uses the equality form of $k$-irreducibility, not the stronger condition involving inclusions in unions, and all images carry the subspace topology. The sliced almost-disjoint family and the seed topology are due to van Douwen; the infinite-chain extension and the analysis of its compact saturated sets provide the kernel counterexample.

Section~\ref{sec:seed} records the seed construction and classifies its proper $k$-irreducible compact subsets. Sections~\ref{sec:extension}--\ref{sec:kernel} define the extension, determine its specialization order, and construct the kernel operator. Sections~\ref{sec:compact}--\ref{sec:image} establish co-sobriety of the extension and failure of co-sobriety of its image.

\section{Preliminaries}

Compactness is used without any separation assumption. We write $\N=\{0,1,2,\ldots\}$ and $[\N]^\omega$ for the family of all infinite subsets of $\N$.

Let $X$ be a topological space. A subset $S\subseteq X$ is \emph{saturated} if it is the intersection of all open sets containing it. If $X$ is $T_0$, this is equivalent to $S$ being an upper set with respect to the specialization order. We write
\[
\sat_X(x)=\{y\in X:x\leq_X y\}
\]
for the saturation of a point $x$.

A nonempty compact saturated set $K\subseteq X$ is called \emph{$k$-irreducible} if
\[
K=K_0\cup K_1,~ K_0,K_1\ \text{compact saturated},
\]
implies $K=K_0$ or $K=K_1$. The family of all such sets is denoted by $\KIrr(X)$. A $T_0$-space $X$ is \emph{co-sober} if for every $K\in\KIrr(X)$ there exists a unique point $x\in X$ such that
\[
K=\sat_X(x).
\]

Recall that a space is a \emph{KC-space} if every compact subset is closed, and it is \emph{anti-Hausdorff} if every two nonempty open subsets intersect. A space is called a \emph{US-space} if every convergent sequence has a unique limit. A \emph{Fr\'echet space} means a Fr\'echet--Urysohn space: for every $A\subseteq X$ and $x\in\cl_X A$, a sequence in $A$ converges to $x$. We shall use the following observation of K\"unzi and van der Zypen \cite[Lemma~1]{KunziZypen2004}.

\begin{lemma}[K\"unzi and van der Zypen]\label{lem:frechet-us-kc}
Every Fr\'echet US-space is a KC-space.
\end{lemma}

\begin{proof}
Let $K$ be a compact subset of a Fr\'echet US-space $X$, and let $x\in\cl_X K$. By the Fr\'echet property there is a sequence $(x_n)$ of points of $K$ converging to $x$. Since $K$ is compact, the sequence has a cluster point $y\in K$. Using the Fr\'echet property again, one can choose a subsequence of $(x_n)$ converging to $y$. Since $X$ is a US-space, limits of convergent sequences are unique, and hence $x=y\in K$. Thus $K$ is closed.
\end{proof}

\begin{lemma}\label{lem:seed-not-cosober}
Let $X$ be a non-singleton compact $T_1$ KC anti-Hausdorff space. Then $X$ is not co-sober.
\end{lemma}

\begin{proof}
Since $X$ is $T_1$, every subset is saturated. Suppose
\[
X=K_0\cup K_1
\]
with $K_0$ and $K_1$ nonempty compact subsets. By the KC property, both $K_i$ are closed. If both are proper, then $X\setminus K_0$ and $X\setminus K_1$ are disjoint nonempty open sets, contradicting anti-Hausdorffness. Hence $X$ is $k$-irreducible. Since $X$ is non-singleton while every point saturation in a $T_1$-space is a singleton, $X$ is not co-sober.
\end{proof}

\section{A van Douwen-type seed space}\label{sec:seed}

We use the sliced maximal almost-disjoint construction of van Douwen. We first state the required combinatorial result and then record the properties of his seed space needed for the extension.

\subsection{A sliced MAD theorem of van Douwen}

The following is the case of van Douwen's slicing theorem for maximal almost-disjoint (MAD) families on a countable set; see \cite[Sections~5--7]{vanDouwen1993}. The underlying set is countable; the family itself is not asserted to be countable.

\begin{theorem}[van Douwen]\label{thm:van-douwen-slicing}
There exist an infinite family $\cA\subseteq[\N]^\omega$ and a map $\ell\colon\cA\to\N$ such that:
\begin{enumerate}
\item[(M1)] distinct $a,b\in\cA$ have finite intersection;
\item[(M2)] for every infinite $S\subseteq\N$ there exists $a\in\cA$ such that $a\cap S$ is infinite;
\item[(M3)] if
\[
\cA(S)=\{a\in\cA:|a\cap S|=\infty\}
\]
is infinite, then
\[
\ell[\cA(S)]=\N.
\]
Equivalently, for every $n\in\N$ there exists $a\in\cA(S)$ with $\ell(a)=n$.
\end{enumerate}
\end{theorem}

\begin{remark}\label{rem:van-douwen-slicing}
Conditions (M1) and (M2) say that $\cA$ is an $\omega$-MAD family. Condition (M3) is precisely the slicing condition denoted by (E) in van Douwen's construction; see \cite[p.~151]{vanDouwen1993}. Its existence in ZFC follows from van Douwen's slicing theorem, in particular Theorem~7.1 of \cite{vanDouwen1993}, after identifying the countably infinite set occurring there with $\N$.
\end{remark}

For $a\in\cA$ and $B\subseteq a$, call
\[
B\cup\{\ell(a)\}
\]
a \emph{branch piece}. Put $Z=\N$ and let the complements of all branch pieces form a subbase for a topology $\tau_Z$ on $Z$.

\begin{proposition}\label{prop:seed-properties}
The space $Z=(\N,\tau_Z)$ is countably infinite, compact, $T_1$, KC and anti-Hausdorff. Moreover every proper closed subset of $Z$ is contained in a set of the form
\[
\bigcup_{a\in F}a\ \cup E,
\]
where $F\subseteq\cA$ and $E\subseteq\N$ are finite.
\end{proposition}

\begin{proof}
We record the parts of the sliced van Douwen construction that will be used later.

First, $\ell$ is onto: apply (M3) to $S=\N$. Consequently each singleton $\{n\}$ occurs as a branch piece, and hence $Z$ is $T_1$.

If $U$ is open and $\ell(a)\in U$, then $a\setminus U$ is finite. It is enough to check this for subbasic open sets. If a subbasic open set containing $\ell(a)$ is the complement of $B\cup\{\ell(b)\}$, then $b\neq a$; by (M1), its complement meets $a$ in only finitely many points. The property is preserved under finite intersections and arbitrary unions. It follows that every branch piece $B\cup\{\ell(a)\}$ is compact: an open cover contains one member covering $\ell(a)$, and that member misses only finitely many points of $B$.

Next, every proper closed subset $C$ is small in the stated sense. Pick $x\notin C$. Since $Z\setminus C$ is open, there is a basic neighbourhood of $x$, given by a finite intersection of subbasic open sets, contained in $Z\setminus C$. Taking complements shows that $C$ is contained in a finite union of branch pieces and hence in $\bigcup_{a\in F}a\cup E$ for suitable finite $F$ and $E$. On the other hand, $Z$ itself is not small: if $F$ and $E$ are finite, choose $a\in\cA\setminus F$ and then a point of $a$ outside $E\cup\bigcup F$, using (M1).

Every proper closed set is compact, because it is a closed subset of a finite union of compact branch pieces and finitely many points. To see that $Z$ is compact, take an open cover and choose one nonempty member $U$. The complement $Z\setminus U$ is a proper closed, hence compact, subset and therefore has a finite subcover; adjoining $U$ gives a finite cover of $Z$.

Van Douwen's argument for the topology determined by such an $\omega$-MAD family and slicing shows that $Z$ is Fr\'echet and US; see \cite[pp.~149--151]{vanDouwen1993}, in particular properties (8)--(10) and condition (E). Hence Lemma~\ref{lem:frechet-us-kc} implies that $Z$ is a KC-space. Finally, if two nonempty open sets $U,V$ were disjoint, then $Z\setminus U$ and $Z\setminus V$ would be proper closed small sets whose union is $Z$, contradicting the fact that $Z$ is not small. Thus $Z$ is anti-Hausdorff.
\end{proof}

\begin{proposition}\label{prop:proper-kirr-seed}
If $K\subsetneq Z$ is a $k$-irreducible compact subset, then $K$ is a singleton.
\end{proposition}

\begin{proof}
By the KC property, $K$ is closed. By Proposition~\ref{prop:seed-properties}, there are a finite set $F\subseteq\cA$ and a finite $E\subseteq Z$ such that
\[
K\subseteq \bigcup_{a\in F}a\cup E.
\]
For each $a\in F$, define
\[
Q_a=
\begin{cases}
K\cap a, & |K\cap a|<\infty,\\[1mm]
(K\cap a)\cup\{\ell(a)\}, & |K\cap a|=\infty.
\end{cases}
\]
If $K\cap a$ is infinite, an injective enumeration of $K\cap a$ converges to $\ell(a)$ by the branch-convergence property in the preceding proof. Since $K$ is closed, $\ell(a)\in K$. Hence every $Q_a$ is contained in $K$, and it is compact: the finite case is immediate and the infinite case is a branch piece. Together with the finitely many singletons belonging to $K\cap E$, these sets give a finite cover of $K$ by compact subsets of $K$. Since $Z$ is $T_1$, all of them are saturated. Repeated use of $k$-irreducibility shows that $K$ must coincide with one of these pieces.

A finite piece containing at least two points splits into two nonempty proper finite compact subsets, so a finite $k$-irreducible piece is a singleton. If $Q_a$ is infinite, choose $x\in K\cap a$ with $x\neq\ell(a)$. Then
\[
Q_a=\{x\}\cup\bigl(((K\cap a)\setminus\{x\})\cup\{\ell(a)\}\bigr)
\]
is a union of two nonempty proper compact subsets, again contradicting the $k$-irreducibility of $K$. Thus only the singleton case is possible.
\end{proof}

By Lemma~\ref{lem:seed-not-cosober} and Proposition~\ref{prop:seed-properties}, $Z$ itself is not co-sober.

\section{The infinite-chain extension}\label{sec:extension}

Fix $z_0\in Z$. Let
\[
Y=Z\sqcup\N.
\]
Write the old points as $\mu(z)$, $z\in Z$, and the new chain points as $p_n$, $n\in\N$.

For $O\subseteq Y$, put
\[
U(O)=\{z\in Z:\mu(z)\in O\},~
A(O)=\{n\in\N:p_n\in O\},
\]
and for $U\subseteq Z$ define
\[
R(U)=Z\setminus(U\cup\{z_0\}).
\]
We declare $O$ to be open if and only if:
\begin{enumerate}
\item[(E1)] $U(O)$ is open in $Z$;
\item[(E2)] $A(O)$ is an upper set in $\N$ with its usual order;
\item[(E3)] if $z_0\in U(O)$, then $A(O)=\N$;
\item[(E4)] if $A(O)\neq\varnothing$, then $U(O)\neq\varnothing$ and $R(U(O))$ is compact in $Z$.
\end{enumerate}

\begin{proposition}\label{prop:extension-topology}
Conditions (E1)--(E4) define a topology on $Y$.
\end{proposition}

\begin{proof}
The empty set and $Y$ satisfy the conditions. Let $O,P$ be open. Conditions (E1)--(E3) are stable under intersections. If $A(O\cap P)\neq\varnothing$, then both $U(O)$ and $U(P)$ are nonempty. Since $Z$ is anti-Hausdorff, $U(O)\cap U(P)\neq\varnothing$, and
\[
R(U(O)\cap U(P))=R(U(O))\cup R(U(P))
\]
is compact. Thus $O\cap P$ is open.

Now let $O=\bigcup_{i\in I}O_i$, where each $O_i$ is open. Conditions (E1)--(E3) are immediate. If $A(O)\neq\varnothing$, choose $i_0$ with $A(O_{i_0})\neq\varnothing$. Then $U(O)\neq\varnothing$ and
\[
R(U(O))=R(U(O_{i_0}))\cap(Z\setminus U(O)).
\]
Since $Z\setminus U(O)$ is closed and $R(U(O_{i_0}))$ is compact, the right-hand side is compact. Hence arbitrary unions are open.
\end{proof}

\begin{lemma}\label{lem:old-embedding}
The map $\mu\colon Z\to Y$ is a topological embedding.
\end{lemma}

\begin{proof}
The old trace $U(O)$ of every open $O\subseteq Y$ is open in $Z$. Conversely, let $U\subseteq Z$ be open. If $z_0\notin U$, then $\mu(U)$ is open in $Y$. If $z_0\in U$, then
\[
\mu(U)\cup\{p_n:n\in\N\}
\]
is open: condition (E4) holds because $R(U)=Z\setminus U$ is a closed subset of the compact space $Z$. In both cases the old trace is exactly $U$.
\end{proof}

\section{Specialization order and the kernel operator}\label{sec:kernel}

\begin{proposition}\label{prop:specialization}
For $z,w\in Z$ and $m,n\in\N$ one has
\begin{align*}
\mu(z)\leq_Y\mu(w) &\iff z=w,\\
p_n\leq_Y p_m &\iff n\leq m,\\
\mu(z_0)\leq_Y p_n &~\text{for every }n,
\end{align*}
and there are no other comparabilities between an old point and a new point.
\end{proposition}

\begin{proof}
The first assertion follows from Lemma~\ref{lem:old-embedding} and the $T_1$ property of $Z$. Condition (E2) gives $p_n\leq_Y p_m$ whenever $n\leq m$. For the converse, the set
\[
O_n=\mu(Z\setminus\{z_0\})\cup\{p_j:j\geq n\}
\]
is open by (E1)--(E4), since its residual old set is empty. If $m<n$, then $p_n\in O_n$ and $p_m\notin O_n$. Thus $p_n\leq_Y p_m$ if and only if $n\leq m$.

Condition (E3) says that every open neighbourhood of $\mu(z_0)$ contains every $p_n$, so $\mu(z_0)\leq_Y p_n$. If $z\neq z_0$, the open set $\mu(Z\setminus\{z_0\})$ contains $\mu(z)$ and no chain point. Conversely, for every old point $\mu(z)$ there is an open set containing a prescribed $p_n$ but omitting $\mu(z)$: use old trace $Z\setminus\{z\}$ when $z\neq z_0$, and $Z\setminus\{z_0\}$ when $z=z_0$, together with a suitable tail. Hence no new point specializes to an old point.
\end{proof}

\begin{corollary}\label{cor:T0}
The space $Y$ is $T_0$.
\end{corollary}

Define
\[
k(\mu(z))=\mu(z)~(z\in Z),~
k(p_n)=\mu(z_0)~(n\in\N).
\]

\begin{proposition}\label{prop:kernel}
The map $k\colon Y\to Y$ is a continuous kernel operator.
\end{proposition}

\begin{proof}
Idempotence is immediate. By Proposition~\ref{prop:specialization},
\[
k(p_n)=\mu(z_0)\leq_Y p_n,
\]
while every old point is fixed. Hence $k(y)\leq_Y y$ for every $y\in Y$.

For continuity, let $O$ be open. If $z_0\notin U(O)$, then
\[
k^{-1}(O)=\mu(U(O)).
\]
If $z_0\in U(O)$, then
\[
k^{-1}(O)=\mu(U(O))\cup\{p_n:n\in\N\}.
\]
Both sets satisfy (E1)--(E4), so $k$ is continuous.
\end{proof}

\section{Compact saturated sets in the extension}\label{sec:compact}

Put
\begin{align*}
D&=\sat_Y(\mu(z_0))
 =\{\mu(z_0)\}\cup\{p_m:m\in\N\},\\
T_n&=\sat_Y(p_n)
 =\{p_m:m\geq n\},\\
P_n&=\{\mu(z):z\neq z_0\}\cup\{p_m:m\geq n\}.
\end{align*}

\begin{lemma}\label{lem:compact-special}
The sets $D$, $T_n$, $P_n$ and $Y$ are nonempty compact saturated subsets of $Y$.
\end{lemma}

\begin{proof}
The sets $D$ and $T_n$ are point saturations and hence supercompact. The set $P_n$ is saturated by Proposition~\ref{prop:specialization}. Let $\mathcal U$ be an open cover of $P_n$ and choose $O_0\in\mathcal U$ with $p_n\in O_0$. By (E2), $O_0$ contains every $p_m$ with $m\geq n$. By (E4), $R(U(O_0))$ is compact. The set $\mu(R(U(O_0)))$ is compact by Lemma~\ref{lem:old-embedding} and is contained in $P_n$. Hence finitely many members of $\mathcal U$ cover it. Together with $O_0$, they cover $P_n$.

Finally, let $\mathcal V$ cover $Y$. Compactness of the old copy $\mu(Z)$ yields finitely many members covering all old points. One of them contains $\mu(z_0)$ and therefore, by (E3), contains every chain point. Thus $Y$ is compact.
\end{proof}

Define
\[
\pi\colon Y\longrightarrow Z,~
\pi(\mu(z))=z,~ \pi(p_n)=z_0.
\]
The same inverse-image computation as in Proposition~\ref{prop:kernel} shows that $\pi$ is continuous.

\begin{lemma}[Full-projection classification]\label{lem:full-projection}
Let $S\subseteq Y$ be saturated and suppose $\pi[S]=Z$. Then either $S=Y$ or $S=P_n$ for some $n\in\N$.
\end{lemma}

\begin{proof}
For each $z\neq z_0$, the fibre $\pi^{-1}(z)$ consists only of $\mu(z)$, so all these old points belong to $S$. If $\mu(z_0)\in S$, saturation and Proposition~\ref{prop:specialization} force every $p_n$ into $S$, and hence $S=Y$.

Assume $\mu(z_0)\notin S$. Since $z_0\in\pi[S]$, the set
\[
N_S=\{m:p_m\in S\}
\]
is nonempty. Saturation makes $N_S$ an upper set of $\N$, so it has a least element $n$ and equals $\{m:m\geq n\}$. Thus $S=P_n$.
\end{proof}

There are strict decompositions
\[
Y=D\cup P_0,~
P_n=T_n\cup P_{n+1}~(n\in\N).
\]
All terms are nonempty compact saturated sets by Lemma~\ref{lem:compact-special}, and both pieces in each union are proper.

\begin{corollary}\label{cor:not-kirr-full}
Neither $Y$ nor any $P_n$ is $k$-irreducible.
\end{corollary}

\section{Co-sobriety of the extension}\label{sec:cosober-extension}

\begin{theorem}\label{thm:Y-cosober}
The space $Y$ is co-sober.
\end{theorem}

\begin{proof}
Let $S\in\KIrr(Y)$. We first show that $\pi[S]$ is a $k$-irreducible compact subset of $Z$. It is nonempty and compact by continuity of $\pi$, and it is saturated because $Z$ is $T_1$.

Suppose
\[
\pi[S]=A\cup B,
\]
where $A,B$ are nonempty compact subsets of $Z$. Since $Z$ is KC, both $A$ and $B$ are closed. Put
\[
S_A=S\cap\pi^{-1}(A),~
S_B=S\cap\pi^{-1}(B).
\]
Both sets are nonempty. They are saturated because $A$ and $B$ are saturated in the $T_1$-space $Z$ and inverse images of saturated sets under continuous maps are saturated. They are compact because they are closed subsets of the compact set $S$. Since
\[
S=S_A\cup S_B,
\]
$k$-irreducibility gives $S=S_A$ or $S=S_B$, and therefore $\pi[S]=A$ or $\pi[S]=B$. Hence $\pi[S]$ is $k$-irreducible.

If $\pi[S]=Z$, Lemma~\ref{lem:full-projection} gives $S=Y$ or $S=P_n$ for some $n$, contradicting Corollary~\ref{cor:not-kirr-full}. Therefore $\pi[S]\neq Z$. By Proposition~\ref{prop:proper-kirr-seed},
\[
\pi[S]=\{q\}
\]
for some $q\in Z$.

If $q\neq z_0$, the fibre $\pi^{-1}(q)$ consists only of $\mu(q)$, and hence
\[
S=\{\mu(q)\}=\sat_Y(\mu(q)).
\]
Suppose $q=z_0$. If $\mu(z_0)\in S$, saturation gives $D\subseteq S$, whereas $S\subseteq\pi^{-1}(z_0)=D$; hence $S=D=\sat_Y(\mu(z_0))$. If $\mu(z_0)\notin S$, then $S$ consists of chain points. Its set of levels is a nonempty upper subset of $\N$, hence equals $\{m:m\geq n\}$ for its least element $n$. Thus
\[
S=T_n=\sat_Y(p_n).
\]
In all cases $S$ is a point saturation. Conversely, every point saturation is supercompact and hence $k$-irreducible: in any cover by compact saturated subsets, a member containing the point contains its entire saturation. Thus the proof gives the exact classification
\[
\KIrr(Y)=\bigl\{\{\mu(q)\}:q\in Z\setminus\{z_0\}\bigr\}
\cup\{D\}\cup\{T_n:n\in\N\}.
\]
By Corollary~\ref{cor:T0}, each representing point is unique, so $Y$ is co-sober.
\end{proof}

\section{The kernel image is not co-sober}\label{sec:image}

By definition,
\[
k(Y)=\mu(Z).
\]
By Lemma~\ref{lem:old-embedding}, the restriction
\[
\mu\colon Z\longrightarrow k(Y)
\]
is a homeomorphism. Since $Z$ is non-singleton, compact, $T_1$, KC and anti-Hausdorff, Lemma~\ref{lem:seed-not-cosober} shows that $Z$, and hence $k(Y)$, is not co-sober.

\begin{theorem}\label{thm:kernel-negative}
There exist a countably infinite compact co-sober $T_0$-space $Y$ and a continuous kernel operator $k\colon Y\to Y$ such that $k(Y)$, endowed with the subspace topology, is not co-sober. Hence the answer to \cite[Question~9.5]{XuZhao2021} is negative.
\end{theorem}

Since $k^2=k$, the image equals the fixed-point set,
\[
k(Y)=\Fix(k),
\]
and $k\colon Y\to k(Y)$ together with the inclusion $k(Y)\hookrightarrow Y$ is a retraction. Moreover the composite inclusion--retraction is exactly the kernel operator $k\leq\mathrm{id}_Y$.

\begin{corollary}\label{cor:retract-negative}
A kernel retract of a co-sober space need not be co-sober. In particular, a retract of a co-sober space need not be co-sober, so the answer to \cite[Question~9.6]{XuZhao2021} is negative.
\end{corollary}

\section{Concluding remarks}

The construction exhibits a sharp contrast between closure and kernel retracts. The image of a co-sober space under a continuous closure operator is co-sober \cite{WenXu2018}, whereas Theorem~\ref{thm:kernel-negative} shows that the analogous assertion for a kernel operator fails. In retract language, closure retracts preserve co-sobriety, while kernel retracts do not.

The infinite chain supplies the decompositions used in this example. Every saturated set whose projection covers all of $Z$ is, by Lemma~\ref{lem:full-projection}, either $Y$ or one of the $P_n$, and each has a nontrivial compact saturated decomposition:
\[
Y=D\cup P_0,~
P_n=T_n\cup P_{n+1}.
\]
All remaining $k$-irreducible compact saturated sets have proper projection to the seed, and the seed classification then collapses that projection to a singleton. These decompositions allow co-sobriety of the ambient space to coexist with non-co-sobriety of its kernel image.

The result suggests a broader study of preservation properties under order-compatible retracts. In particular, it is natural to compare closure retracts, kernel retracts and arbitrary retracts for other non-Hausdorff properties arising in domain theory.

\section*{Declaration of competing interest}
The authors declare that they have no known competing financial interests or personal relationships that could have appeared to influence the work reported in this paper.

\section*{Data availability}
No research data were used in this theoretical study.

\section*{Declaration of generative AI and AI-assisted technologies in the manuscript preparation process}
During the preparation of this manuscript, the authors used ChatGPT (OpenAI) for language editing, organization of the exposition, literature cross-checking, and assistance with checking mathematical arguments. After using this tool, the authors reviewed and edited the content as needed and take full responsibility for the content of the publication.

\end{document}